\documentstyle[12pt,leqno]{article}
 \newcommand{\beq}{\begin{equation}}
\newcommand{\eeq}{\end{equation}}
 \newcommand{\beth}{\begin{theo}}
\newcommand{\eth}{\end{theo}}

\newcommand{\Subset}{\subset\subset}

\newcommand{\taups}{\tau(\Psi)}
\newcommand{\Psux}{{ \Psi_{u,x}}}
\newcommand{\Psuo}{{ \Psi_{u,0}}}
\newcommand{\psux}{{ \psi_{u,x}}}
\newcommand{\taux}{\tau(u,x)}
\newcommand{\nux}{\nu(u,x)}
\newcommand{\nuxa}{\nu(u,x,a)}
\newcommand{\taupsux}{\tau(\Psi_{u,x})}
\newcommand{\taupsuxo}{\tau(\Psi_{u,x},0)}
\newcommand{\tmu}{T_m u}
\newcommand{\cL}{{Log}}
\newcommand{\cE}{{Exp}}

\newcommand{\okn}{1\le k\le n}

\newcommand{\Znp}{{ \bf Z}_+^n}
\newcommand{\Rn}{{ \bf R}^n}
\newcommand{\Rnm}{{ \bf R}_-^n}
\newcommand{\Rnp}{{ \bf R}_+^n}
\newcommand{\Cn}{{ \bf C}^n}
\newcommand{\Cp}{{ \bf C}^p}

\newtheorem{theo}{Theorem}
\newtheorem{prop}{Proposition}
\begin{document}

 \title{Newton numbers and residual measures of plurisubharmonic functions}
\author{\sc Alexander RASHKOVSKII  }
\date{}
\maketitle

  {\small {\bf Abstract.}
We study the masses charged by $(dd^cu)^n$ at isolated singularity
points of plurisubharmonic functions $u$. It is done by means of the
local indicators of plurisubharmonic functions introduced in
\cite{keyLeR}. As a consequence,
bounds for the masses are obtained in terms
of the directional Lelong numbers of $u$, and
the notion of the Newton number for a holomorphic mapping is
extended to arbitrary plurisubharmonic functions.
We also describe the local
indicator of $u$ as the logarithmic tangent to $u$.


}

\vskip 0.5cm

\section{Introduction}

The principal information on local behaviour of a subharmonic function
$u$ in the complex plane can be obtained by studying its
Riesz measure $\mu_u$.
If $u$ has a logarithmic singularity at a point $x$,
the main term of its asymptotics near $x$ is
$\mu_u(\{x\})\log|z-x|$.
For plurisubharmonic functions $u$ in $\Cn$, $n>1$,
the situation is not so simple. The local propoperties of $u$ are
controlled by the current $dd^cu$
(we use the notation $d=\partial + \bar\partial,\ d^c= ( \partial
-\bar\partial)/2\pi i$) which cannot charge isolated points.
The trace measure $\sigma_u=dd^cu\wedge\beta_{n-1}$
of this current is precisely the Riesz measure of $u$; here
$\beta_p=(p!)^{-1}2^{p}(dd^c|z|^2)^p$ is the volume element of $\Cp$.
A significant role is played by the Lelong numbers $\nux$
of the function $u$ at  points $x$:
$$
\nux=\lim_{r\to 0}\,(\tau_{2n-2}r^{2n-2})^{-1}\sigma_u[B^{2n}(x,r)],
$$
where $\tau_{2p}$ is the volume of the unit ball $B^{2p}(0,1)$ of $\Cp$.
If $\nux>0$ then $\nux\log|z-x|$ gives an upper bound for $u(z)$ near
$x$, however the difference between these two functions can be
comparable to $\log|z-x|$.

Another important object generated by the current $dd^cu$
is the Monge-Amp\'ere measure $(dd^cu)^n$.
For the definition and basic facts on the complex Monge-Amp\'ere
operator $(dd^c)^n$ and Lelong numbers,
we refer the reader to the books
\cite{keyL2}, \cite{keyLeG} and \cite{keyKl},
and for more advanced results, to \cite{keyD1}. Here we mention that
$(dd^cu)^n$ cannot be defined for all plurisubharmonic functions $u$,
however if $u\in PSH(\Omega)\cap L_{loc}^\infty(\Omega\setminus K)$
with $K\Subset\Omega$, then $(dd^cu)^n$ is well defined as a positive
closed current of the bidimension $(0,0)$ (or, which is the same, as
a positive measure) on $\Omega$. This measure cannot charge
pluripolar subsets of $\Omega\setminus K$, and it can have positive
masses at points of $K$, e.g. $(dd^c\log|z|)^n=\delta(0)$, the Dirac
measure at $0$, $|z|=(\sum|z_j|^2)^{1/2}$. More generally, if
$f:\Omega\to{\bf C}^N,\ N\ge n$, is a holomorphic mapping with
isolated zeros at $x^{(k)}\in\Omega$ of multiplicities $m_k$, then
$(dd^c\log|f|)^n|_{x^{(k)}}=m_k\,\delta(x^{(k)})$. So, the masses
of $(dd^cu)^n$ at isolated points of singularity of $u$ (the residual
measures of $u$) are of especial importance.

Let a plurisubharmonic function $u$ belong to $L_{loc}^\infty
(\Omega\setminus\{x\})$;
its residual mass at the point $x$ will be denoted by $\taux$:
$$
\taux=(dd^cu)^n|_{\{x\}}.
$$
The problem under consideration is evaluation of this value.

The following well-known relation compares $\taux$ with the Lelong
number $\nux$:
\beq
\taux\ge [\nux]^n.
\label{eq:3}
\eeq
The equality in (\ref{eq:3}) means that, roughly speaking, the
function $u(z)$ behaves near $x$ as $\nux\log|z-x|$.
In many cases however relation (\ref{eq:3}) is not optimal; e.g.
for
\beq
u(z)=\sup\{\log|z_1|^{k_1},\log|z_2|^{k_2}\},\quad k_1>k_2,
\label{eq:4}
\eeq
$\tau(u,0)=k_1k_2>k_2^2=[\nu(u,0)]^2$.

As follows from the Comparison Theorem due to Demailly (see Theorem~A
below), the residual mass is determined by asymptotic behaviour
of the function near its singularity, so one needs to find
appropriate characteristics for the behaviour. To this end,
a notion of local indicator was proposed
in \cite{keyLeR}. Note that $\nux$ can be calculated as $$
\nux=\lim_{r\to -\infty}r^{-1}\sup\{v(z):\: |z-x|\le e^r\}
= \lim_{r\to -\infty}r^{-1}{\cal M}(u,x,r),
$$
where ${\cal M}(u,x,r)$ is the mean value of $u$ over the sphere
$|z-x|=e^r$, see \cite{keyKis1}. In \cite{keyKis2}, the {\it refined},
or {\it directional, Lelong numbers} were introduced as
\begin{eqnarray}
\nuxa &=& \lim_{r\to -\infty}r^{-1}\sup\{v(z):\: |z_k-x_k|\le
           e^{ra_k},\ \okn\}\nonumber\\
    &=& \lim_{r\to -\infty}r^{-1}g(u,x,ra),
\label{eq:dir}
\end{eqnarray}
where $a=(a_1,\ldots,a_n)\in\Rnp$ and $g(u,x,b)$ is the mean value
of $u$ over the set $\{z:\: |z_k-x_k|=\exp{b_k},\ \okn\}$.
For $x$ fixed, the collection $\{\nuxa\}_{a\in\Rnp}$ gives a more
detailed information about the function $u$ near $x$ than $\nux$
does, so one can expect for a more precise bound for $\taux$ in terms
of the directional Lelong numbers.
It was noticed already in \cite{keyKis2} that $a\mapsto\nuxa$ is a
concave function on $\Rnp$. In \cite{keyLeR}, it was observed that
this function produces the following plurisubharmonic function
$\Psux$ in the unit polydisk $D=\{y\in\Cn: |y_k|<1,\ \okn\}$:
$$
\Psux(y)=-\nu(u,x,(-\log|y_k|)),
$$
the {\it local indicator} of the function $u$ at $x$. It is the
largest negative plurisubharmonic function in $D$ whose directional
Lelong numbers
at $0$ coincide with those of $u$ at $x$,
$(dd^c\Psux)^n=\taupsuxo\,\delta(0)$, and finally,
\beq
\taux\ge\taupsuxo,
\label{eq:5}
\eeq
so the singularity of $u$ at $x$ is controlled by its indicator
$\Psux$.

Since $\taupsuxo\ge[\nu(\Psux,0)]^n=[\nux]^n$, (\ref{eq:5}) is
a refinement of (\ref{eq:3}). For the function $u$ defined by
(\ref{eq:4}), $\tau(\Psi_{u,0},0)=k_1k_2=\tau(u,0)>[\nu(u,0)]^2$.

Being a function of a quite simple nature, the indicator can produce
effective bounds for residual measures of plurisubharmonic functions.
In Theorems~1--3 of the present paper we study the values
$N(u,x):=\taupsuxo$, {\it the Newton numbers of $u$ at} $x$; the
reason for this name is explained below.  We obtain,  in particular,
the following bound for $\taux$ (Theorem~\ref{theo:4}):
$$
\taux\ge\frac{[\nuxa]^n}{a_1\ldots a_n}\quad \forall a\in\Rnp;
$$
it reduces to (\ref{eq:3}) when $a_1=\ldots=a_n=1$.
For $n$ plurisubharmonic functions $u_1,\ldots,u_n$
in general position (see the definition below), we estimate the
measure $dd^c\Psi_{u_1,x}\wedge\ldots\wedge dd^c\Psi_{u_n,x}$ and
prove the similar relation (Theorem~\ref{theo:6})
\beq
dd^cu_1\wedge\ldots\wedge dd^cu_n|_{\{x\}}\ge
\frac{\prod_j\nu(u_j,x,a)}{a_1\ldots a_n}\quad\forall a\in \Rnp.
\label{eq:7}
\eeq

The main tool used to obtain these bounds is the Comparison
Theorem due to Demailly. To formulate it we give the following

{\it Definition 1.} A $q$-tuple of plurisubharmonic functions
$u_1,\ldots,u_q$
is said to be {\it in general position} if their unboundedness
loci $A_1,\ldots,A_q$ satisfy the following condition:
for all choices of indices $j_1<\ldots<j_k,\ k\le q$,
the $(2q-2k+1)$-dimensional Hausdorff measure of $A_{j_1}\cap
\ldots\cap A_{j_k}$ equals zero.

\medskip
{\bf Theorem A} (Comparison Theorem, \cite{keyD1}, Th. 5.9).
{\sl Let $n$-tuples of plurisubharmonic functions $u_1,\ldots,u_n$
and $v_1,\ldots,
v_n$ be in general position on a neighbourhood of a point $x\in\Cn$.
Suppose that $u_j(x)=-\infty,\ 1\le j\le n$, and
$$
\limsup_{z\to x}\frac{v_j(z)}{u_j(z)}=l_j<\infty.
$$
Then}
$$
dd^cv_1\wedge\ldots\wedge dd^cv_n|_{\{x\}}\le
l_1\ldots l_n \,dd^cu_1\wedge\ldots\wedge dd^cu_n|_{\{x\}}.
$$

\medskip
We also obtain a geometric
interpretation for the value $N(u,x)$ (Theorem~\ref{theo:7}). Let
$\Theta_{u,x}$ be the set of points $b\in\overline{\Rnp}$
such that $\nuxa\ge\langle b,a\rangle$ for some $a\in\Rnp$,
then
\beq
\taux\ge N(u,x)=n!\,Vol(\Theta_{u,x}).
\label{eq:8}
\eeq
In many cases the folume of $\Theta_{u,x}$ can be easily calculated, so
(\ref{eq:8}) gives an effective formula for $N(u,x)$.

To illustrate these results, consider functions
$u=\log|f|,\ f=(f_1,\ldots,f_n)$ being an equidimensional
holomorphic mapping with an isolated zero at a point $x$.
It is probably the only class of functions whose residual measures
were studied in details before.
In this case, $\taux$ equals $m$,
the multiplicity of $f$ at $x$, and
\beq
\nu(\log|f|,x,a)=I(f,x,a):=\inf\{\langle a, p \rangle:\:
p\in\omega_x\}
\label{eq:index}
\eeq
where
%\beq
$$
\omega_x=\{p\in\Znp:\: \sum_j \left| \frac{\partial^p f_j}{\partial
        z^p} (x)\right|\neq 0\}
$$
(see \cite{keyLe}). For polynomials $F:\Cn\to{\bf C}$, the value
$I(F,x,a)$ is a known object ({\it the index of $F$ at $x$ with
respect to the weight} $a$) used in number theory (see e.g.
\cite{keyLa}).

Relation (\ref{eq:5}) gives us $m=\tau(\log|f|,x)\ge N(\log|f|,x)$.
In general, the
value $N(\log|f|,x)$ is not comparable to $m_1\ldots m_n$ with $m_j$ the
multiplicity of the function $f_j$: for $f(z)=(z_1^2+z_2, z_2)$ and
$x=0$, $m_1m_2=1<2=N(\log|f|,x)=m$ while for $f(z)=(z_1^2+z_2, z_2^3)$,
$N(\log|f|,x)=2<3=m_1m_2<6=m$.
A more sharp bound for $m$ can be obtained by (\ref{eq:7})
with $u_j=\log|f_j|,\ 1\le j\le n$.
In this case, the left-hand side of (\ref{eq:7}) equals $m$,
and its right-hand side with $a_1=\ldots=a_n$ equals $m_1\ldots m_n$.
For the both above examples of the mapping $f$, the supremum of the
right-hand side of (\ref{eq:7}) over $a\in\Rnp$ equals $m$. For
$a_1,\ldots,a_n$ rational, relation (\ref{eq:7}) is a known bound
for $m$ via the multiplicities of weighted homogeneous initial
Taylor polynomials of $f_j$ with respect to the weights
$(a_1,\ldots,a_n)$ (\cite{keyAYu}, Th. 22.7).

Recall that the convex hull $\Gamma_+(f,x)$ of the set
$\bigcup_p \{p+\Rnp\}$, ${p\in\omega_x}$
is called {\it the Newton polyhedron} of $(f_1,\ldots,f_n)$ at $x$,
the union $\Gamma(f,x)$ of the compact faces of the boundary
of $\Gamma_+(f,x)$ is called {\it the Newton boundary}
of $(f_1,\ldots,f_n)$ at $x$, and the value
$N_{f,x}=n!\,Vol(\Gamma_-(f,x))$ with $\Gamma_-(f,x)=
\{\lambda t: t\in\Gamma(f,x),\ 0\le
\lambda\le 1\}$
is called
{\it the Newton number} of $(f_1,\ldots,f_n)$ at $x$
(see \cite{keyKu1}, \cite{keyAYu}). The relation
\beq
m\ge N_{f,x}
\label{eq:9}
\eeq
was established by
A.G.~Kouchnirenko \cite{keyKu} (see also \cite{keyAYu}, Th. 22.8).
Since $\Theta_{\log|f|,x}=
\Gamma_-(f,x)$, (\ref{eq:9}) is a particular case of the relation
(\ref{eq:8}). It is the reason to
call $N(u,x)$ {the Newton number} of $u$ at $x$.

These observations show that the technique of plurisubharmonic
functions (and local indicators in particular) is quite a powerful
tool to produce, in a unified and simple way, sharp bounds for the
multiplicities of holomorphic mappings.

Finally, we obtain a description for the indicator $\Psux(z)$ as
the weak limit of the functions $m^{-1}u(x_1+z_1^m,\ldots,
x_n+z_n^m)$ as $m\to\infty$ (Theorem~\ref{theo:indctr}), so $\Psux$
can be viewed as the tangent (in the logarithmic coordinates) for the
function $u$ at $x$.  Using this approach we obtain a sufficient
condition, in terms of ${\cal C}_{n-1}$-capacity, for the residual
mass $\taux$ to coincide with the Newton number of $u$ at $x$
(Theorem~\ref{theo:converg}).

%%%%%%%%%%%%%%%%%%%%%%%%%%%%%%%%%%%%

\section{Indicators and their masses}

We will use the following notations. For a domain $\Omega$ of $\Cn$,
$PSH(\Omega)$ will denote the class of all plurisubharmonic functions
on $\Omega$, $PSH_-(\Omega)$ the subclass of the nonpositive
functions, and $PSH(\Omega,x)=PSH(\Omega)\cap L_{loc}^\infty(\Omega
\setminus\{x\})$ with $x\in\Omega$.

Let $D=\{z\in\Cn:\: |z_k|<1,\ \okn\}$ be the unit polydisk,
$D^*=\{z\in D:\: z_1\cdot\ldots\cdot z_n\neq 0\}$,
${\bf R}_\pm^n=\{t\in\Rn:\: \pm t_k>0\}$. By $CNVI_-(\Rnm)$ we denote
the collection of all nonpositive convex functions on $\Rnm$
increasing in each variable $t_k$.
The mapping $\cL:D^*\to\Rnm$ is defined as $\cL(z)=(\log|z_1|,\ldots,
\log|z_n|)$, and $\cE:\Rnm\to D^*$ is given by $\cE(t)=(\exp
t_1,\ldots, \exp t_n)$.

A function $u$ on $D^*$ is called $n$-{\it circled} if
\beq
u(z)=u(|z_1|,\ldots,
|z_n|),
\label{eq:circ}
\eeq
i.e. if $\cL^*\cE^*u=u$. Any $n$-circled function
$u\in PSH_-(D^*)$ has a unique extension to the whole polydisk $D$
keeping the property (\ref{eq:circ}). The class of such functions
will be denoted by $PSH_-^c(D)$. The cones $CNVI_-(\Rnm)$ and
$PSH_-^c(D)$ are isomorphic: $u\in PSH_-^c(D) \iff \cE^*u\in
CNVI_-(\Rnm)$, $h\in CNVI_-(\Rnm) \iff \cL^*h\in PSH_-^c(D)$.

\medskip
{\it Definition 2} \cite{keyLeR}. A function $\Psi\in PSH_-^c(D)$ is
called an {\it indicator} if its convex image $\cE^*\Psi$ satisfies
\beq
\cE^*\Psi(ct)=c\,\cE^*\Psi(t)\quad \forall c>0,\ \forall t\in\Rnm.
\label{eq:odn}
\eeq

\medskip
The collection of all indicators will be denoted by $I$. It is a
convex subcone of $PSH_-^c(D)$, closed in ${\cal D}'$ (or
equivalently, in $L_{loc}^1(D)$). Besides, if $\Psi_1,\Psi_2
\in I$ then $\sup\{\Psi_1,\Psi_2\}\in I$, too.

Every indicator is locally bounded in $D^*$. In what follows we will
often consider indicators locally bounded in $D\setminus\{0\}$; the
class of such indicators will be denoted by $I_0$: $I_0=I\cap
PSH(D,0)$.

An example of indicators can be given by the functions
$$
\varphi_a(z)=\sup_k\,a_k\log|z_k|,\ a_k\ge 0.
$$
If all $a_k>0$,
then $\varphi_a\in I_0$.

\begin{prop}
Let $\Psi\in I_0,\ \Psi\not\equiv 0$. Then
\begin{enumerate}
\item[(a)] there exist $\nu_1,\ldots,\nu_n>0$ such that
\beq
\Psi(z)\ge \varphi_\nu(z)\quad\forall z\in D;
\label{eq:lbound}
\eeq
\item[(b)] $\Psi\in C(\overline D\setminus\{0\}),\ \Psi|_{\partial
D}=0$;
\item[(c)] the directional Lelong numbers $\nu(\Psi,0,a)$ of $\Psi$
at the origin with respect to $a\in\Rnp$ (\ref{eq:dir})
are
\beq
\nu(\Psi,0,a)=-\Psi(\cE(-a)),
\label{eq:dirpsi}
\eeq
and its Lelong number $\nu(\Psi,0)=-\Psi(e^{-1},\ldots,e^{-1})$;
\item[(d)] $(dd^c\Psi)^n=0$ on $D\setminus\{0\}$.
\end{enumerate}
\label{prop:1}
\end{prop}

{\it Proof.} Let $\Psi_k(z_k)$ denote the restriction of
the indicator $\Psi(z)$
to the disk $D^{(k)}=\{z\in D: z_j=0\ \forall j\neq k\}$. By
monotonicity of $\cE^*\Psi$, $\Psi(z)\ge\Psi_k
(z_k)$. Since $\Psi_k$
is a nonzero indicator in the disk $D^{(k)}\subset {\bf C}$,
$\Psi_k(z_k)=\nu_k\log|z_k|$ with some $\nu_k>0$, and (a) follows.

As $\cE^*\Psi\in C(\Rnm),\ \Psi\in C(D^*)$. Its continuity in
$D\setminus\{0\}$ can be shown by induction in $n$. For $n=1$ it is
obvious, so assuming it for $n\le l$, consider any point $z^0\neq 0$
with $z_j^0=0$. Let $z^{s}\to z^0$, then the points $\tilde z^s$
with $\tilde z_j^s=0$ and $\tilde z_m^s=z_m^s,\ m\neq j$,
also tend to $z^0$, and by the induction hypothesis,
$\Psi(\tilde z^s)\to\Psi(\tilde z^0)=\Psi(z^0)$. So,
$\liminf_{s\to\infty}\Psi(z^s)\ge\lim_{s\to\infty}\Psi(\tilde z^s)
=\Psi(z^0)$, i.e. $\Psi$ is lower semicontinuous and hence continuous
at $z^0$. Continuity of $\Psi$ up to $\partial D$ and the boundary
condition follow from (\ref{eq:lbound}).

Equality (\ref{eq:dirpsi}) is an immediate consequence of the
definition of the directional Lelong numbers (\ref{eq:dir}) and
the homogeneity condition (\ref{eq:odn}). The relation
$\nux=\nu(u,x,(1,\ldots,1))$ \cite{keyKis2} gives us the desired
expression for $\nu(\Psi,0)$.

Finally, statement (d) follows from the homogeneity condition
(\ref{eq:odn}), see \cite{keyLeR}, Proposition 4.

\medskip
For functions $\Psi\in I_0$, the complex Monge-Amp\'ere operator
$(dd^c\Psi)^n$ is well defined and gives a nonnegative measure on
$D$. By Proposition \ref{prop:1},
$$
(dd^c\Psi)^n=\taups\delta(0)
$$
with some constant $\taups\ge 0$ which is strictly positive unless
$\Psi\equiv 0$. In this section, we will study the value $\taups$.

An upper bound for $\taups$ is given by

\begin{prop}
For $\Psi\in I_0$,
\beq
\taups\le\nu_1\ldots\nu_n
\label{eq:lbound1}
\eeq
with  $ \nu_1,\ldots,\nu_n$ the same as in Proposition \ref{prop:1},
(a).
\label{prop:2}
\end{prop}

{\it Proof.}
The function $\varphi_\nu(z)\in I_0$, and
(\ref{eq:lbound}) implies
$$
\limsup_{z\to 0}\frac{\Psi(z)}{\varphi_\nu(z)}\le 1,
$$
so (\ref{eq:lbound1}) follows by Theorem A.

\medskip
To obtain a lower bound for $\taups$, we need
 a relation between
$\Psi(z)$ and $\Psi(z^0)$ for $z,z^0\in D$. Denote
$$
\Phi(z,z^0)=\sup_k\,\frac{\log|z_k|}{|\log|z_k^0||},\quad z\in D,
\ z^0\in D^*.
$$
Being considered as a function of $z$ with $z^0$ fixed,
$\Phi(z,z^0)\in I_0$.

\begin{prop}
For any $\Psi\in I$, $\Psi(z)\le |\Psi(z^0)|\Phi(z,z^0)\quad
\forall z\in D,\ z^0\in D^*$.
\label{prop:3}
\end{prop}

{\it Proof.} For a fixed $z^0\in D^*$ and $t^0={Log}(z^0)$, define
$u=|\Psi(z^0)|^{-1}{Exp}^*\Psi$ and $ v={Exp}^*\Phi=\sup_k\,
t_k/|t_k^0|$.
It suffices to establish the inequality $u(t)\le v(t)$ for all
$t\in\Rnm$ with $t_k^0<t_k<0$, $1\le k\le n$.
Given such a $t$, denote $\lambda_0=[1+v(t)]^{-1}$. Since
$\{t^0+\lambda(t-t^0):\: 0\le\lambda\le\lambda_0\}\subset
\overline{\Rnm}$, the functions
$u_t(\lambda):=u(t^0+\lambda(t-t^0))$ and
$v_t(\lambda):=v(t^0+\lambda(t-t^0))$ are well defined on
$[0,\lambda_0]$. Furthermore, $u_t$ is convex and $v_t$ is linear
there, $u_t(0)=v_t(0)=-1$, $u_t(\lambda_0)\le v_t(\lambda_0)=0$.
It implies $u_t(\lambda)\le v_t(\lambda)$ for all
$\lambda\in [0,\lambda_0]$. In particular, as $\lambda_0>1$,
$u(t)=u_t(1)\le v_t(1)=v(t)$, that completes the proof.

\medskip
Consider now the function
$$
P(z)=-\prod_{\okn}|\log|z_k||^{1/n}\in I.
$$

\beth
For any $\Psi\in I_0$,
\beq
 \taups\ge \left|\frac{\Psi(z^0)}{P(z^0)}\right|^n \quad\forall
z^0\in D^*.
\label{eq:ubound1}
\eeq
\label{theo:1}
\eth

{\it Proof.} By Proposition \ref{prop:3},
$$
\frac{\Psi(z)}{\Phi(z,z^0)}\le |\Psi(z^0)|\quad \forall z\in D,\
z^0\in D^*.
$$
By Theorem A,
$$
(dd^c\Psi)^n\le |\Psi(z^0)|^n(dd^c\Phi(z,z^0))^n,
$$
and the statement follows from the fact that
$$
(dd^c\Phi(z,z^0))^n= \prod_{\okn}|\log|z_k^0||^{-1}=|P(z^0)|^{-n}.
$$

\medskip
{\it Remarks.}
1. One can consider the value
\beq
A_\Psi=\sup_{z\in D} \left|\frac{\Psi(z)}{P(z)}\right|^n;
\label{eq:apsi}
\eeq
by Theorem \ref{theo:1},
\beq
\taups\ge A_\Psi.
\label{eq:bound}
\eeq

2. Let $I_{0,M}=\{\Psi\in I_0:\: \taups\le M\},\ M>0$. Then
(\ref{eq:ubound1}) gives the lower bound for the class $I_{0,M}$:
$$
\Psi(z)\ge M^{1/n}P(z)\quad\forall z\in D,\ \forall \Psi\in I_{0,M}.
$$

\medskip
Let now $\Psi_1,\ldots,\Psi_n\in I$ be in general position in the
sense of Definition 1. Then the current $\bigwedge_kdd^c\Psi_k$ is
well defined, as well as $(dd^c\Psi)^n$ with $\Psi=\sup_k\,\Psi_k$.
Moreover, we have

\begin{prop}
If $\Psi_1,\ldots,\Psi_n\in I$ are in general position, then
\beq
\bigwedge_kdd^c\Psi_k=0\ {\rm on\ }D\setminus\{0\}.
\label{eq:prod}
\eeq
\label{prop:nwedge}
\end{prop}

{\it Proof.}
For $\Psi_1,\ldots,\Psi_n\in I_0$, the statement
follows from Proposition
\ref{prop:1}, (d), and the polarization formula
\beq
\bigwedge_kdd^c\Psi_k=\frac{(-1)^n}{n!}\sum_{j=1}^n(-1)^j
\sum_{1\le i_1<\ldots<i_j\le n}\left(dd^c\sum_{k=1}^j \Psi_{j_k}
\right)^n.
\label{eq:polar}
\eeq
When the only condition on $\{\Psi_k\}$ is to be in general position,
we can replace $\Psi_k(z)$ with $\Psi_{k,N}(z)=\sup\{\Psi_k(z),
N\sup_j \log|z_j|\}\in I_0$ for which
$\bigwedge_kdd^c\Psi_{k,N}=0$ on $D\setminus\{0\}$. Since
$\Psi_{k,N}\searrow\Psi_{k}$ as $N\to\infty$, it gives us
(\ref{eq:prod}).

\medskip
The mass of $\bigwedge_kdd^c\Psi_k$ will be denoted by $\tau(\Psi_1,
\ldots,\Psi_n)$.

\beth
Let $\Psi_1,\ldots,\Psi_n\in I$ be in general position,
$\Psi=\sup_k\,\Psi_k$. Then
\begin{enumerate}
\item[(a)] $\taups\le \tau(\Psi_1,\ldots,\Psi_n)$;
\item[(b)] $\tau(\Psi_1,\ldots,\Psi_n)\ge |P(z^0)|^{-n}
\prod_k |\Psi_k(z^0)|\quad\forall z^0\in D^*$.
\end{enumerate}
\label{theo:2}
\eth

{\it Proof.} Since
$$
\frac{\Psi(z)}{\Psi_k(z)}\le 1\ \forall z\neq 0,
$$
statement (a) follows from Theorem A.

Statement (b) results from Proposition \ref{prop:3} exactly like the
statement of Theorem~\ref{theo:1}.

%%%%%%%%%%%%%%%%

\section{Geometric interpretation}

In this section we study the masses $\taups$ of indicators $\Psi\in
I_0$ by means of their convex images $\cE^*\Psi\in CNVI_-(\Rnm)$.

Let $V\in PSH_-^c(rD)\cap C^2(rD),\ r<1$, and $v=\cE^*V\in CNVI_-
\left(({\bf R}_-+\log r)^n\right)$. Since
$$
\frac{\partial^2 V(z)}{\partial z_j\partial \bar z_k}=
\left.{1\over 4z_j\bar z_k}\frac{\partial^2 v(t)}{\partial t_j\partial
t_k}\right|_{t=\cL(z)},\quad z\in rD^*,
$$
$$
\det\left(\frac{\partial^2 V(z)}{\partial z_j\partial \bar z_k}
\right)=
4^{-n}|z_1\ldots z_n|^{-2}
\det\left.\left(\frac{\partial^2 v(t)}{\partial t_j\partial
t_k}\right)\right|_{t=\cL(z)}.
$$

By setting $z_j=\exp\{t_j+i\theta_j\},\ 0\le\theta\le 2\pi$, we get
$\beta_n(z)= |z_1\ldots z_n|^{2}dt\,d\theta$, so
\beq
(dd^cV)^n=n!\left({2\over\pi}\right)^n
\det\left(\frac{\partial^2 V}{\partial z_j\partial \bar z_k}
\right)\,\beta_n=
n!\,(2\pi)^{-n}
\det\left(\frac{\partial^2 v}{\partial t_j\partial
t_k}\right)\,dt\,d\theta.
\label{eq:real}
\eeq

Every function $U\in PSH_-^c(D)\cap L^\infty(D)$ is the limit of
a decreasing sequence of functions $U_l\in PSH_-^c(E)\cap C^2(E)$
on an $n$-circled domain $E\Subset D$, and by the convergence
theorem for the complex Monge-Amp\'ere operators,
\beq (dd^cU_l)^n|_{E}\longrightarrow (dd^cU)^n|_{E}.
\label{eq:conv1}
\eeq
On the other hand, for $u_l=\cE^*U_l$ and $u=\cE^*U$,
\beq
\det\left.\left(\frac{\partial^2 u_l}{\partial t_j\partial
t_k}\right)\,dt\right|_{\cL(D^*\cap E)}
\longrightarrow \left.{\cal MA}[u]\right|_{\cL(D^*\cap E)},
\label{eq:conv2}
\eeq
the {\it real} Monge-Amp\'ere operator of $u$ \cite{keyRaT}.

Since  $(dd^cU_l)^n$ and $(dd^cU)^n$ cannot charge pluripolar sets,
(\ref{eq:real}) with $V=U_l$ and (\ref{eq:conv1}), (\ref{eq:conv2})
imply
$$
(dd^cU)^n(E)= n!\,(2\pi)^{-n}{\cal MA}[u]\,d\theta\, (\cL(E)\times
[0,2\pi]^n)
$$
for any $n$-circled Borel set $E\in D$, i.e.
\beq
(dd^cU)^n(E)= n!\,{\cal MA}[u](\cL(E)).
\label{eq:tr}
\eeq
This relation allows us to calculate $\taups$ by using the technique
of real Monge-Amp\'ere operators in $\Rn$ (see \cite{keyRaT}).

Let $\Psi\in I$. Consider the set
$$
B_\Psi=\{a\in\Rnp:\: \langle a,t\rangle\le \cE^*\Psi(t)\quad
\forall t\in\Rnm\}
$$
and define
$$
\Theta_\Psi=\overline{\Rnp\setminus B_\Psi}.
$$
Clearly, the set $B_\Psi$ is convex,
so $\cE^*\Psi$ is the restriction of its support function to $\Rnm$.
 If $\Psi\in I_0$, the set $\Theta_\Psi $ is bounded.
Indeed, $a\in\Theta_\Psi$ if and only if $\langle a,t^0\rangle\ge
\cE^*\Psi(t^0)$ for some $t^0\in\Rnm$, that implies
$|a_j|\le |\cE^*\Psi(t^0)/t_j^0|\ \forall j$. By Proposition
\ref{prop:1}, (a), $|\cE^*\Psi(t^0)|\le\nu_j|t_j|$ and therefore
$|a_j|\le\nu_j\ \forall j$.

Given a set $F\in\Rn$, we denote its Eucledean volume by $Vol(F)$.

\beth
$\forall \Psi\in I_0$,
\beq
\taups=n!\,Vol(\Theta_\Psi).
\label{eq:vol}
\eeq
\label{theo:3}
\eth

{\it Proof.} Denote $U(z)=\sup\,\{\Psi(z),-1\}\in PSH_-^c(D)\cap
C(D)$, $u=\cE^*U\in CNVI_-(\Rnm)$. Since $U(z)=\Psi(z)$
near $\partial D$,
$$
\taups=\int_D (dd^cU)^n.
$$
Furthermore, as $ (dd^cU)^n=0$ outside the set $E=\{z\in D:\:
\Psi(z)=-1\}$,
\beq
\taups=\int_E (dd^cU)^n.
\label{eq:u1}
\eeq
In view of (\ref{eq:tr}),
\beq
\int_E (dd^cU)^n=n!\,\int_{\cL(E)} {\cal MA}[u].
\label{eq:u2}
\eeq
As was shown in \cite{keyRaT}, for any convex function $v$ in a
domain $\Omega\subset\Rn$,
\beq
\int_{F} {\cal MA}[v]=Vol(\omega(F,v))\quad\forall F\subset\Omega,
\label{eq:u3}
\eeq
where
$$
\omega(F,v)=\bigcup_{t^0\in F} \{a\in\Rn:\: v(t)\ge v(t^0)+
\langle a,t-t^0\rangle\quad\forall t\in\Omega\}
$$
is the gradient image of the set $F$ for the surface $\{y=v(x),\
x\in\Omega\}$.

We claim that
\beq
\Theta_\Psi=\omega(\cL(E),u).
\label{eq:u4}
\eeq

Observe that
$$
\Theta_\Psi=\{a\in\overline{\Rnp}:\: \sup_{\psi(t)=-1} \langle
a,t\rangle\ge -1\}
$$
where $\psi=\cE^*\Psi$.

If $a\in \omega(\cL(E),u)$, then for some $t^0\in\Rnm$ with
$\psi(t^0)=1$ we have $\langle a,t^0\rangle\ge\langle a,t\rangle$
for all $t\in\Rnm$ such that $\psi(t)<-1$. Taking here
$t_j\to -\infty$ we get $a_j\ge 0$, i.e. $a\in\overline{\Rnp}$.
Besides, $\langle a,t^0\rangle\ge\langle a,t\rangle -1-\psi(t)$
for all $t\in\Rnm$ with $\psi(t)>-1$, and applying this for $t\to 0$
we derive $\langle a,t^0\rangle\ge -1$. Therefore, $a\in\Theta_\Psi$
and $\Theta_\Psi\supset\omega(\cL(E),u)$.

Now we prove the converse inclusion. If  $a\in\Theta_\Psi\cap\Rnp$,
then
$$
\sup\{\langle a,t^0\rangle:\: t^0\in \cL(E)\}\ge -1.
$$
Let $t$ be such that $\psi(t)=-\delta>-1$, then $t/\delta\in \cL(E)$
and thus
\begin{eqnarray*}
\langle a,t\rangle -1-\psi(t) &=& \delta\langle a,t/\delta\rangle
    -1+\delta \le
  \delta\sup_{t^0\in \cL(E)}\langle a,t^0\rangle -1+\delta\\
   & \le & \sup_{t^0\in \cL(E)}\langle a,t^0\rangle
 = \sup_{z^0\in E}\langle a,\cL(z^0)\rangle.
\end{eqnarray*}
Since $E$ is compact, the latter supremum is attained at some point
$\hat z^0$. Furthermore, $\hat z^0\in E\cap D^*$ because $a_k\neq 0,
\ \okn$. Hence  $\sup_{t^0\in \cL(E)}\langle a,t^0\rangle=
\langle a,\hat t^0\rangle$ with $\hat t^0=\cL(z^0)\in\Rnm$, so that
$a\in\omega(\cL(E),u) $ and
$\Theta_\Psi\cap\Rnp\subset\omega(\cL(E),u)$.
Since $\omega(\cL(E),u)$ is closed, this implies
$\Theta_\Psi=\omega(\cL(E),u)$, and (\ref{eq:u4}) follows.

Now relation (\ref{eq:vol}) is a consequence of
(\ref{eq:u1})--(\ref{eq:u4}). The theorem is proved.

\medskip
Note that the value $\tau(\Psi_1,\ldots,\Psi_n)$ also can be
expressed in geometric terms. Namely, if $\Psi_1,\ldots,\Psi_n
\in I_0$, the polarization formula (\ref{eq:polar}) gives us,
by Theorem \ref{theo:3},
$$
\tau(\Psi_1,\ldots,\Psi_n)= (-1)^n\sum_{j=1}^n(-1)^j
\sum_{1\le i_1<\ldots<i_j\le n} Vol(\Theta_{\sum_k\Psi_{j_k}}).
$$

We can also give an interpretation for the bound (\ref{eq:bound}).
Write $A_\Psi$ from (\ref{eq:apsi}) as
\beq
A_\Psi =\sup_{a\in\Rnp}\frac{|\psi(-a)|^n}{a_1\ldots a_n}
= \sup_{a\in\Rnp}\left|\psi(-{a/ a_1})\ldots
\psi(-{a/ a_n})\right|,
\label{eq:apsi1}
\eeq
$\psi=\cE^*\Psi$. For any $a\in\Rnp$, the point $a^{(j)}$ whose
$j$th coordinate equals $|\psi(-{a/ a_j})|$ and the others are
zero, has the property  $\langle a^{(j)},-a\rangle=\psi(-a)$. This
remains true for every convex combination $\sum\rho_j a^{(j)}$ of the
points $a^{(j)}$, and thus $r\sum\rho_j a^{(j)}\in\Theta_\Psi$
with any $r\in [0,1]$.
Since

$ (n!)^{-1} \left|\psi(-{a/ a_1})\ldots\psi(-{a/ a_n})\right|$
is the volume of the simplex generated by the points $0,a^{(1)},
\ldots,a^{(n)}$, we see from (\ref{eq:apsi1}) that $(n!)^{-1}A_\Psi$
is the supremum of the volumes of all simplices contained in
$\Theta_\Psi$.

Besides, $(n!)^{-1}[\nu(\Psi,0)]^n$ is the volume of the simplex
$$
\{a\in\overline{\Rnp}:\: \langle a,(1,\ldots,1)\rangle \le
\nu(\Psi,0)\}\subset\Theta_\Psi.
$$
It is a geometric description for
the "standard" bound $\taups\ge [\nu(\Psi,0)]^n$.

%%%%%%%%%%%%%%%%%%%%%

\section{Singularities of plurisubharmonic functions}
Let $u$ be a plurisubharmonic function in a domain $\Omega\subset\Cn$,
and $\nuxa$
be its directional Lelong number (\ref{eq:dir}) at $x\in\Omega$
with respect to $a\in\Rnp$. Fix a point $x$. As is known
\cite{keyKis2}, the function $a\mapsto\nuxa$ is a concave function
on $\Rnp$. So, the function
$$
\psux(t):=-\nu(u,x,-t),\quad t\in\Rnm,
$$
belongs to $CNVI_-(\Rnm)$ and thus
$$
\Psux:=\cL^*\psux\in PSH_-^c(D).
$$
Moreover, due to the positive homogeneity of $\nuxa$ in $a$,
$\Psux\in I$. The function $\Psux$ was introduced in \cite{keyLeR}
as {\it (local) indicator} of $u$ at $x$. According to
(\ref{eq:dir}),
\begin{eqnarray*}
\Psux(z) &=&\lim_{R\to +\infty}R^{-1}\sup\{u(y):\: |y_k-x_k|\le
               |z_k|^R,\ \okn\}\\
  &=& \lim_{R\to +\infty}R^{-1} {1\over (2\pi)^n}
       \int_{[0,2\pi]^n}
       u(x_k+|z_k|^R e^{i\theta_k})\,d\theta_1 \ldots d\theta_n.
\end{eqnarray*}
Clearly, $\Psux\equiv 0$ if and only if $\nu(u,x)=0$.
It is easy to see that $\Psi(\Phi,0)=\Phi\ \forall\Phi\in I$.
In particular,
\beq
\nuxa=\nu(\Psux,0,a)=-\Psux(\cE(-a))\quad \forall a\in\Rnp.
\label{eq:rel}
\eeq
So, the results of the previous sections can be applied to study
the  directional Lelong numbers of arbitrary plurisubharmonic
functions.

\begin{prop} {\rm (cf. \cite{keyKis3}, Pr. 5.3)}\
For any $u\in PSH(\Omega)$,
$$
\nuxa\ge \nu(u,x,b)\sup_k\,{a_k\over b_k}\quad\forall x\in\Omega,
\ \forall a,b\in \Rnp.
$$
\label{prop:4}
\end{prop}

{\it
 Proof.} In view of (\ref{eq:rel}), the relation follows from
Proposition \ref{prop:3}.

\medskip
Given $r\in\Rnp$ and $z\in\Cn$, we denote $r^{-1}=(r_1^{-1},\ldots,
r_n^{-1})$ and $r\cdot z=(r_1z_1,\ldots,r_nz_n)$.

\begin{prop} {\rm(\cite{keyLeR})}.
If $u\in PSH(\Omega)$
then
\beq
u(z)\le\Psux(r^{-1}\cdot z)+\sup\, \{u(y):\: y\in
D_r(x)\}
\label{eq:A}
\eeq
for all $z\in D_r(x)= \{y:\:|y_k-x_k|\le
r_k,\ \okn\}\Subset\Omega$.
\label{prop:5}
\end{prop}

{\it Proof.}
Let us assume for simplicity $x=0$, $D_r(0)=D_r$.
Consider the function $v(z)=u(r\cdot z) -\sup \{u(y):\: y\in D_r\}
\in PSH_-(D)$.
The function $g_v(R,t):= \sup\{v(z):\: |z_k|\le\exp\{Rt_k\},\
\okn\}$ is convex in $R>0$ and $t\in\Rnm$, so for $R\to\infty$
\beq
\frac{g_v(R,t)-g_v(R_1,t)}{R-R_1}\nearrow\psi_{v,0}(t),
\label{eq:mon}
\eeq
$\psi_{v,0}=\cE^*\Psi_{v,0}$.

For $R=1,\ R_1\to 0$, (\ref{eq:mon}) gives us $g_v(1,t)\le
\psi_{v,0}(t)$ and thus (\ref{eq:A}).
The proposition is proved.

\medskip
Let $\Omega_k(x)$ be the connected component of the set
$\Omega\cap\{z\in\Cn:\: z_j=x_j\quad\forall j\neq k\}$
containing the point $x$. If for some $x\in\Omega$,
$u|_{\Omega_k(x)}\not\equiv -\infty\quad\forall k$,
then $\Psux\in I_0$. For example, it is the case for $u\in
PSH(\Omega,x)$.

If $u\in PSH(\Omega,x)$, the measure $(dd^cu)^n$ is defined on
$\Omega$. Its residual mass at $x$ will be denoted by $\taux$:
$$
\taux= (dd^cu)^n|_{\{x\}}.
$$
Besides, the idicator $\Psux\in I_0$. Denote $N(u,x)=\taupsux$.

\begin{prop} {\rm(\cite{keyLeR}, Th. 1)}.
If  $u\in PSH(\Omega,x)$, then
$\taux\ge N(u,x)$.
\label{prop:6}
\end{prop}

{\it Proof.} Inequality (\ref{eq:A}) implies
$$
\limsup_{z\to x}\frac{\Psux(r^{-1}\cdot(z-x))}{u(z)}\le 1,
$$
and since
$$
\lim_{y\to 0} \frac{\Psux(r^{-1}\cdot y))}{\Psux(y)}=1\quad
\forall r\in\Rnp,
$$
the statement follows from Theorem A.

\medskip
So, to estimate $\taux$ we may apply the bounds for $\taupsux$
from the previous section.

\beth
If  $u\in PSH(\Omega,x)$, then
$$
\taux\ge\frac{[\nuxa]^n}{a_1\ldots a_n}\quad \forall a\in\Rnp;
$$
in other words,
$\taux\ge A_{u,x}$
where $A_{u,x}=A_{\Psux}$ is defined by (\ref{eq:apsi}).
\label{theo:4}
\eth

{\it Proof.} The result follows from Theorem \ref{theo:1} and
Proposition \ref{prop:6}.

\medskip
Let now $u_1,\ldots,u_n\in PSH(\Omega)$ be in general position in the
sense of Definition 1. Then the current $\bigwedge_k dd^cu_k$ is
defined on $\Omega$ (\cite{keyD1}, Th. 2.5); denote its residual mass
at a point $x$ by $\tau(u_1,\ldots,u_n;x)$.
Besides, the $n$-tuple
of the indicators $\Psi_{u_k,x}$ is in general position, too, that
implies $\bigwedge_k dd^c\Psi_{u_k,x} =\tau(\Psi_{u_1,x},\ldots,
\Psi_{u_n,x})\,\delta(0)$ (Proposition \ref{prop:nwedge}).

In view of Theorem A and Proposition \ref{prop:5} we have

\beth
$\tau(u_1,\ldots,u_n;x) \ge \tau(\Psi_{u_1,x},\ldots,
\Psi_{u_n,x})$.
\label{theo:5}
\eth

Now Theorems \ref{theo:2} and \ref{theo:5} give us

\beth
\beq
\tau(u_1,\ldots,u_n;x) \ge
\frac{\prod_j\nu(u_j,x,a)}{a_1\ldots a_n}\quad\forall a\in \Rnp.
\label{eq:C}
\eeq
\label{theo:6}
\eth

{\it Remark.} For $a_1=\ldots=a_n=1$, inequality (\ref{eq:C}) is
proved in \cite{keyD1}, Cor. 5.10.

\medskip
By combination of Proposition \ref{prop:6} and Theorem \ref{theo:3}
we get

\beth
For $u\in PSH(\Omega,x)$,
\beq
\taux\ge N(u,x)=n!\,V(\Theta_{u,x})
\label{eq:D}
\eeq
with
$$
\Theta_{u,x}=\{b\in\Rnp:\: \sup_{\sum a_k=1}[\nuxa-\langle b,a
\rangle]\ge 0\}.
$$
\label{theo:7}
\eth

\medskip
{\it Remark on holomorphic mappings.} Let $f=(f_1,\ldots,f_n)$ be a
holomorphic mapping of a neighbourhood $\Omega$ of the origin into
$\Cn$, $f(0)=0$ be its isolated zero. Then in a subdomain
$\Omega'\subset\Omega$
the zero sets $A_j$ of the functions $f_j$ satisfy the
conditions
$$
A_1\cap\ldots\cap A_n\cap\Omega'=\{0\},\quad
{\rm codim}\: A_{j_1}\cap\ldots\cap A_{j_k}\cap\Omega'\ge k
$$
for all choices
of indices $j_1<\ldots<j_k,\ k\le n$. Denote $u=\log|f|,\
u_j=\log|f_j|$. Then, as is known,
$\tau(u,0)=\tau(u_1,\ldots,u_n;0)=m_f$, the multiplicity of $f$ at
$0$.  For $a=(1,\ldots,1)$, $\nu(u_j,0,a)$ equals $m_j$, the
multiplicity of the function $f_j$ at $0$. Therefore, (\ref{eq:C})
with $a=(1,\ldots,1)$ gives us the standard bound $m_f\ge m_1\ldots
m_n$.

For $a_j$ rational, (\ref{eq:C}) is the known estimate of $m_f$ via
the multiplicities of weighted homogeneous initial Taylor polynomials
for $f_j$ (see e.g. \cite{keyAYu}, Th. 22.7). Indeed, due to the
positive homogeneity of the directional Lelong numbers, we can take
$a_j\in\Znp$. Then by (\ref{eq:index}), $\nu(u_j,0,a)$ is equal to
the multiplicity of the function $f_j^{(a)}(z)=f_j(z^a)$.

We would also like to mention that (\ref{eq:C})
gives a lower bound for the Milnor number
$\mu(F,0)$ of a singular point $0$ of a holomorphic function $F$
(i.e. for the multiplicity of the isolated zero of the mapping
$f={\rm grad}\, F$ at $0$) in terms of the indices $I(F,0,a)$
(\ref{eq:index}) of
$F$.  Since $I({\partial F/ \partial z_k},0,a)\ge I(F,0,a)-a_k$,
$$
\mu(F,0)\ge\prod_{\okn}\left( {I(F,0,a)\over a_k}-1\right).
$$

Finally, as follows from (\ref{eq:index}), the set
$\overline{\Rnp\setminus\Theta_{u,0}}$ is the Newton
polyhedron for the system $(f_1,\ldots,f_n)$ at $0$ (see
Introduction).
Therefore, $n!\,V(\Theta_{u,0})$ is the Newton number of
$(f_1,\ldots,f_n)$ at
 $0$, and (\ref{eq:D}) becomes the bound
for $m_f$ due to A.G.~Kouchnirenko (see \cite{keyAYu}, Th. 22.8).
So, for any plurisubharmonic function $u$, we will call the value
$N(u,x)$ {\it the Newton number of $u$ at $x$}.

%%%%%%%%%%%%%%%%%%%

\section{Indicators as logarithmic tangents}

Let $u\in PSH(\Omega,0),\ u(0)=-\infty$. We will consider the
following problem: under what conditions on $u$,
its residual measure equals its Newton number?

Of course, the relation
\beq
\exists\lim_{z\to 0} \frac{u(z)}{\Psuo(z)}=1
\label{eq:regul}
\eeq
is sufficient, however it seems to be too restrictive.
On the other hand, as the example $u(z)=\log(|z_1+z_2|^2+|z_2|^4)$
shows, the condition
$$
\lim_{\lambda\to 0} \frac{u(\lambda z)}{\Psuo(\lambda z)}=1
\quad \forall z\in\Cn\setminus\{0\}
$$
does not guarantee %(\ref{eq:rav}).
the equality $\tau(u,0)=N(u,0)$.

To weaken (\ref{eq:regul}) we first give another description for the
local indicators.
In \cite{keyKis4}, a compact family of plurisubharmonic functions
$$
u_r(z)=u(rz)-\sup\{u(y):\: |y|<r\}_{r>0}
$$
was considered
and the limit sets, as $r\to 0$, of such families were described.
In particular, the limit set need not consist of a single function,
so a plurisubharmonic function can have several (and thus infinitely
many) tangents.
Here we consider another family generated by a plurisubharmonic
function $u$.

Given $m\in{\bf N}$ and $z\in\Cn$, denote
$z^m=(z_1^m,\ldots, z_n^m)$ and set
$$
\tmu(z)=m^{-1}u(z^m).
$$
Clearly, $T_mu\in PSH(\Omega\cap D)$ and $T_mu\in
PSH_-(\overline{D_r})$ for any $r\in \Rnp\cap D^*$ (i.e. $0<r_k<1$)
for all $m\ge m_0(r)$.

\begin{prop}
The family $\{T_mu\}_{m\ge m_0(r)}$ is compact in $L_{loc}^1(D_r)$.
\label{prop:comp}
\end{prop}

{\it Proof.} Let $M(v,\rho)$ denote the mean value of
a function $v$ over the set
$\{z:\: |z_k|=\rho_k,\ \okn\}$, $0<\rho_k\le r_k$, then
$M(T_mu,\rho) =m^{-1} M(u,\rho^m)$.
The relation
\beq
m^{-1} M(u,\rho^m)\nearrow \Psuo(\rho)\ {\rm as}\ m\to\infty
\label{eq:c1}
\eeq
implies $M(T_mu,\rho)\ge M(T_{m_0}u,\rho)$.
Since $T_mu\le 0$ in $D_r$, it proves the compactness.

\beth
\begin{enumerate}
\item[(a)]
$T_mu\to\Psuo$ in $L_{loc}^1(D)$;
\item[(b)]
if $u\in PSH(\Omega,0)$ then
$(dd^c T_mu)^n\to\tau(u,0)\,\delta(0)$.
\label{theo:indctr}
\end{enumerate}
\eth

{\it Proof.} Let $g$ be a partial limit of the sequence $T_mu$,
that is $T_{m_s}u\to g$ as $s\to\infty$ for some sequence $m_s$.
For the function $v(z)=\sup\,\{u(y):\: |y_k|\le |z_k|,\ \okn\}$
and any $r\in\Rnp\cap D^*$ we have by (\ref{eq:A})
$$
\tmu(z)\le \left(T_mv\right)(z)\le\Psuo(r^{-1}\cdot z)
$$
and thus
\beq
g(z)\le\Psuo(z)\ \forall z\in D.
\label{eq:c2}
\eeq

On the other hand, the convergence of $T_{m_s}u$ to $g$ in $L^1$
implies $M(T_{m_s}u,r)\to M(g,r)$ (\cite{keyHo}, Prop. 4.1.10).
By (\ref{eq:c1}), $M(T_{m_s}u,r)\to\Psuo(r)$, so
$M(g,r)=\Psuo(r)$ for every $r\in\Rnp\cap D^*$. Being compared with
(\ref{eq:c2}) it gives us $g\equiv\Psuo$, and the statement {\it (a)}
follows.

To prove {\it (b)} we observe that for each $\alpha\in (0,1)$
$$
\int_{\alpha D} (dd^c T_mu)^n= \int_{\alpha^m D} (dd^c u)^n
\to \tau(u,0)
$$
as $m\to\infty$, and for $0<\alpha<\beta<1$
$$
\lim_{m\to\infty} \int_{\beta D\setminus\alpha D} (dd^c T_mu)^n
= \lim_{m\to\infty}\left[\int_{\beta^m D} (dd^cu)^n
-\int_{\alpha^m D} (dd^cu)^n\right]=0.
$$
The theorem is proved.

\medskip

So, Theorem \ref{theo:indctr} shows us that $\tau(u,0)=N(u,0)$
if an only if  $(dd^c T_mu)^n\to (dd^c \Psuo)^n$. And now we are
going to find conditions for this convergence.

Recall  the definition of the inner ${\cal C}_{n-1}$-capacity
introduced in \cite{keyX}:
for any Borel subset $E$ of a domain $\omega$,
$$
{\cal C}_{n-1}(E,\omega)=\sup\,\{\int_E(dd^cv)^{n-1}\wedge\beta_1:\:
v\in PSH(\omega),\ 0<v<1\}.
$$
It was shown in \cite{keyX} that convergence of uniformly bounded
plurisubharmonic functions $v_j$ to $v$ in ${\cal C}_{n-1}$-capacity
imples $(dd^cv_j)^n\to (dd^cv)^n$. In our situation, neither $T_mu$
 nor $\Psuo$ are bounded, so we will modify the construction from
 \cite{keyX}.

 Set
 $$
 E(u,m,\delta)=\{z\in D\setminus\{0\}:\:
 \frac{\tmu (z)}{\Psuo(z)}>1+\delta\}, \quad m\in{\bf N},\ \delta>0.
 $$

 \beth
 Let $u\in PSH(\Omega,0), \rho\in(0,1/4), N>0$, and a sequence
 $m_s\in{\bf N}$ be such that
 \begin{enumerate}
 \item[1)]
 $u(z)>-Nm_s$ on a neighbourhood of the sphere $\partial
 B_{\rho^{m_s}}$, $\forall s$;
 \item[2)]
 $\lim_{s\to\infty}{\cal C}_{n-1}(B_\rho\cap E(u,m_s,\delta), D)=0
 \quad\forall\delta>0$.
 \end{enumerate}
 Then $(dd^c T_mu)^n\to (dd^c\Psuo)^n$ on $D$.
 \label{theo:converg}
 \eth

 {\it Proof.}
 Without loss of generality we can take $u\in PSH_-(D,0)$.
Consider the functions
$v_s(z)=\max\,\{T_{m_s}u(z),-N\}$ and
 $v=\max\,\{\Psuo(z),-N\}$. We have $v_s=T_{m_s}u$ and $v=\Psuo$
 on a neighbourhood of $\partial B_\rho$, $v_s=v=-N$ on a
 neighbourhood of $0$, $v_s\le v$ on $B_\rho$, and
 $v_s\ge (1+\delta)v$ on $B_\rho\setminus E(u,m_s,\delta)$.

 We will prove the relations
 \beq
 (dd^c v_s)^k\wedge (dd^c v)^l\to (dd^c v)^{k+l}
 \label{eq:induct}
 \eeq
 for $k=1,\ldots,n,\ l=0,\ldots,n-k$.
 As a consequence, it will give us the statement of the theorem.
 Indeed, by  Theorem \ref{theo:indctr},
 $$
 \int_{B_\rho} (dd^cv_s)^n= \int_{B_\rho} (dd^cT_{m_s}u)^n
 \to \tau(u,0)
 $$
 while
$$
 \int_{B_\rho} (dd^cv)^n= \int_{B_\rho} (dd^c\Psuo)^n =N(u,0),
 $$
 and (\ref{eq:induct}) with $k=n$ provides the coincidence of the
 right-hand sides of these relations and thus the convergence
 of $(dd^c T_mu)^n$ to $(dd^c\Psuo)^n$.

 We prove (\ref{eq:induct}) by induction in $k$. Let $k=1,
 0\le l\le n-1, \delta>0$. For any test form $\phi\in
 {\cal D}_{n-l-1,n-l-1}(B_\rho)$,
$$
 \left|\int dd^cv_s\wedge(dd^cv)^l\wedge\phi
   -\int(dd^cv)^{l+1}\wedge\phi\right|
   = \left|\int (v-v_s)(dd^cv)^l\wedge dd^c\phi\right|
$$
\begin{eqnarray*}
 \le & C_\phi & \int_{B_\rho} (v-v_s)(dd^cv)^l\wedge\beta_{n-l} \\
   = & C_\phi & \left[\int_{B_\rho\setminus E_{s,\delta}}
                + \int_{B_\rho\cap E_{s,\delta}} \right]
                 (v-v_s)(dd^cv)^l\wedge\beta_{n-l} \\
   = & C_\phi & [I_1(s,\delta)+I_2(s,\delta)],
 \end{eqnarray*}
 where, for brevity, $E_{s,\delta}=E(u,m_s,\delta)$.

 We have
 $$
 I_1(s,\delta)\le\delta \int_{B_\rho}
|v|(dd^cv)^l\wedge\beta_{n-l} \le C\delta
 $$
 with a constant $C$ independent of $s$, and
 \begin{eqnarray*}
 I_2(s,\delta) &\le & N\int_{B_\rho\cap E_{s,\delta}}
 (dd^cv)^l\wedge\beta_{n-l}\\
      &\le & C(N,\rho,l)\cdot {\cal C}_{n-1}(B_\rho\cap E_{s,\delta}, D)
 \longrightarrow 0.
 \end{eqnarray*}
 Since $\delta>0$ is arbitrary, it proves (\ref{eq:induct}) for
 $k=1$.

 Let us now have got (\ref{eq:induct}) for $k=j$ and $0\le l\le n-j$.
 For $\phi\in {\cal D}_{n-l-j,n-l}(B_\rho)$,
 $$
 \int  (dd^cv_s)^{j+1}\wedge(dd^cv)^l\wedge\phi=
 \int  (dd^cv_s)^{j}\wedge(dd^cv)^{l+1}\wedge\phi
 $$
 $$
 + \int \left[ (dd^cv_s)^{j+1}\wedge(dd^cv)^l-
 (dd^cv_s)^{j}\wedge(dd^cv)^{l+1}\right]\wedge\phi.
 $$
 The first integral in the right-hand side converges to
 $\int(dd^cv)^{l+j+1}\wedge\phi$ by the induction assumption.
 The second integral can be estimated similarly to the case $k=1$:
 $$
 \left|\int \left[ (dd^cv_s)^{j+1}\wedge(dd^cv)^l-
 (dd^cv_s)^{j}\wedge(dd^cv)^{l+1}\right]\wedge\phi\right|
 $$
\begin{eqnarray*}
\le & C_\phi & \left[\int_{B_\rho\setminus E_{s,\delta}}
               + \int_{B_\rho\cap E_{s,\delta}} \right]
               (v-v_s)(dd^cv_s)^j(dd^cv)^l\wedge\beta_{n-j-l} \\
 =  & C_\phi & [I_3(s,\delta)+I_4(s,\delta)].
 \end{eqnarray*}
Since  $(dd^cv_s)^{j}\wedge(dd^cv)^l\to(dd^cv)^{j+l}$,
  $$
 \int (dd^cv_s)^j(dd^cv)^l\wedge\beta_{n-j-l} \le C\quad\forall s
 $$
 and
 $$
 I_3(s,\delta)\le \delta
 \int_{B_\rho} |v|(dd^cv_s)^j(dd^cv)^l\wedge\beta_{n-j-l}
\le CN\delta.
 $$
 Similarly,
 \begin{eqnarray*}
 I_4(s,\delta) &\le & N
   \int_{B_\rho\cap E_{s,\delta}} (dd^cv_s)^j(dd^cv)^l
       \wedge\beta_{n-j-l} \\
  &\le & C(N,\rho,j,l)\cdot {\cal C}_{n-1}(B_\rho\cap E_{s,\delta}, D)
 \longrightarrow 0,
 \end{eqnarray*}
 and (\ref{eq:induct}) is proved.

\vskip 0.5cm

Mathematical Division, Institute for Low Temperature Physics
\par
47 Lenin Ave., Kharkov 310164,
Ukraine

\vskip0.1cm

E-mail: \quad rashkovskii@ilt.kharkov.ua,
 \quad rashkovs@ilt.kharkov.ua

\end{document}